\documentclass[10pt]{article}
\usepackage{mathrsfs}
\usepackage{amsmath}
\usepackage{amssymb}
\usepackage{graphicx}
\newtheorem{theorem}{\scshape \mdseries  Theorem}[section]
\newtheorem{lemma}[theorem]{\scshape \mdseries  Lemma}
\newtheorem{coro}[theorem]{\scshape \mdseries  Corollary}

\newtheorem{rem}[theorem]{\scshape \mdseries  Remark}

\begin{document}

\title{\sf On the reciprocal degree distance of graphs with cut vertices or edges  \thanks{
      Supported by Natural Science Research Foundation on Department of Education of Anhui Province(KJ2013A196). }}
\author{Xiao-Xin Li \thanks{E-mail address:
 lxx@czu.edu.com }
\\
  {\small  \it Department of Mathematics, Chizhou University, Chizhou 247000, P.R. China}
  }
\date{}
\maketitle

\noindent {\bf Abstract}\ \ As an additive weight version of the Harary index, the reciprocal degree distance of a simple connected graph $G$ is defined as $RDD(G)=\sum_{\{u,v\}\subseteq V(G)}\frac{d_G(u)+d_G(v)}{d_G(u,v)},$ where $d_G(u)$ is the degree of $u$ and $d_G(u,v)$ is the distance between $u$ and $v$ in $G$. In this paper, the extremal graph with the maximum RRD-value among all the graphs of order $n$ with given cut vertices or edges is characterized. In addition, an upper bounds on the reciprocal degree distance in terms of the number of cut edges is provided.



\section{Introduction}
Chemical graphs are models of molecules in which atoms and chemical bonds are represented by vertices and edges of a graph, respectively. Chemical graph theory is a branch of mathematical chemistry concerning the study of chemical graph. A graph invariant (also known as molecular descriptor or topological index) is a function on a graph that does not depend on a labeling of its vertices. The chemical information derived through topological index has been found useful in chemical documentation, isomer discrimination, structure property correlations, etc \cite{ash}. Hundreds of graph invariants of molecular graphs are studied in chemical graph theory. Many of them are based on the graph distance, see \cite{xu2} and references therein, another large group is based on the vertex degree, see \cite{fur} and references therein. In addition, several graph invariants are based on both the vertex degree and the graph distance, see \cite{xu3} and the references therein. In this paper, we are interested in a distance-degree-based graph invariant which is called the reciprocal degree distance of a graph.

Let $G$ be a simple connected graph with vertex set $V(G)$ and edge set $E(G)$.
$d_G(v)$ denotes the degree of a vertex $v$ in $G$ and $d_G(u,v)$ denotes the distance between two vertices $u$ and $v$ in $G$.

For a connected graph $G$, one of the oldest and well-known distance-based graph invariants is Wiener index, denoted by $W(G)$, which is introduced by Wiener \cite{wie} in 1947 and defined as the sum of distance over all unordered vertex pairs in $G$, i.e., $$W(G)=\sum_{\{u,v\}\subseteq V(G)} d_G(u,v).$$

Another distance-based graph invariant is Harary index, denoted by $H(G)$, which is defined as the sum of reciprocals of distances between all pairs of vertices in $G$, i.e. $$H(G)=\sum_{\{u,v\}\subseteq V(G)} \frac{1}{d_G(u,v)}.$$

In 1994, a degree-weighted version of Wiener index called degree distance or Schultz molecular toplogical index was proposed by Dobrynin and Kochetova \cite{dob} and Gutman \cite{gut} independently, which is defined for a connected graph $G$ as $$DD(G)=\sum_{\{u,v\}\subseteq V(G)}(d_G(u)+d_G(v))d_G(u,v).$$

The interested readers may consult \cite{dob2,gut2,gut3,hua} for Wiener index, \cite{das,diu,feng,he} for Harary index and \cite{buc,dan,ili2,tom1,tom3,tom4,tom5} for degree distance.

Similarly, a degree-weighted version of Harary index called reciprocal degree distance was proposed by Alizadeh et al. \cite{ali} in 2013 and Hua and Zhang \cite{hua2} in 2012 independently, which is defined for a connected graph $G$ as $$RDD(G)=\sum_{\{u,v\}\subseteq V(G)}\frac{d_G(u)+d_G(v)}{d_G(u,v)}.$$
It was shown in \cite{ali} that this index can be used as an efficient measuring tool in the study of complex networks.

In general, for a given graph $G$, $RDD(G)$ is not always easily calculated. So it makes sense to determine the bounds of $RDD(G)$ or to characterize the graphs with extremal reciprocal degree distance among a given class of graphs. In \cite{hua2}, Hua and Zhang established various lower and upper bounds for the reciprocal degree distance among various given class of graphs including tree, unicyclic graph, cactus and given pendent vertices, independence number, chromatic number, vertex connectivity and edge connectivity. Li and Meng \cite{li1} characterized the extremal graphs among $n$ vertex trees with given some graphic parameters such as pendants, matching number, domination number, diameter, vertex bipartition, and determined some sharp upper bounds of trees. Li et al. \cite{li2} determined the maximum RDD-value among all the graphs of diameter $d$ and the connected bipartite graphs with given matching number (resp. vertex connectivity). However, to our best knowledge, the RDD-value of connected
graphs with cut vertices or cut edges has not been considered by other authors so
far. Motivated by the above results, we proceed with the study on the reciprocal degree distance. In this paper, we characterize the unique graph with the maximum RDD-value among all graphs with a given number of cut vertices or edges, and provide an upper bound of the reciprocal degree distance in terms of the number of cut edges.

\section{Preliminaries}

As usual, we begin with some notations and terminology. Let $G$ be a graph, $N_G(v)$ denotes the neighborhood of $v$ in $G$, so $|N_G(v)|=d_G(v)$. A vertex $v$ of $G$ is called pendent if $d_G(v)=1$, and the edge incident with $v$ is called a pendent edge of $G$. A pendent path at $v$ of $G$ is a path in which no vertex other than $v$ is incident with any edge of $G$ outside the path, where the degree of $v$ is at least three. A cut vertex (edge, respectively) of a graph is a vertex (an edge, respectively) whose removal increases the number of components of the graph. A block of a connected graph is defined to be a maximum connected subgraph without cut vertices. A block containing only one cut vertex is called a pendent block, and a block containing only an unique vertex is called trivial. Denote by $P_s=Pv_1v_2...v_s$ a path on vertices $v_1,v_2,...,v_s$ with edges $v_iv_{i+1}$ for $i=1,2,...,s-1$, and denote by $K_n$ a complete graph with order $n$. For simplicity, we denote by $\mathcal{G}_{n,k}$($\overline{\mathcal{G}}_{n,k}$, respectively) the set of connected graphs of order $n$ with $k$ cut vertices(edges, respectively), and denote by $G_{n,k}$ the graph obtained from the complete graph $K_{n-k}$ by adding $n-k$ paths of almost equal lengths to its vertices respectively, denote by $\overline{G}_{n,k}$ the graph obtained from the complete graph $K_{n-k}$ by attaching $k$ pendent vertices to one vertex.

For a subset $V_1\subset V(G)$, let $G-V_1$ be the subgraph of $G$ obtained by deleting the vertices of $V_1$ together with the edges incident with them. If $V_1=\{v\}$, we denote by $G-v$ for simplicity. Similarly, for a subset $E_1\subset E(G)$, let $G-E_1$ be the subgraph of $G$ obtained by deleting the edges of $E_1$. For a subset $E_2\subset E(\overline{G})$, let $G+E_2$ be the graph obtained from $G$ by adding the edges of $E_2$, where $\overline{G}$ is the complement of $G$.  If $E_1=\{e\}$ ($E_2=\{e\}$, respectively), we denote by $G-e$ ($G+e$, respectively) for simplicity.

Note that in any disconnected graph $G$, the distance of any two vertices from two distinct components is infinite.
Therefore its reciprocal can be viewed as 0. Thus, we can define validly the reciprocals degree distance of disconnected graph $G$ as follows:
$$RDD(G)=\sum_{i=1}^k RDD(G_i),$$ where $G_1, G_2, \ldots , G_k$ are the components of $G$.

Let $D_G(u)=\sum\limits_{v \in V(G)\backslash \{u\}}\frac{1}{d_G(u,v)}.$ By direct calculation or from \cite{hua2}, we can get $$RDD(G)=\sum_{u \in V(G)}d_G(u)D_G(u).$$

Let $e=(u,v)$ be an edge of $G$. The removal of $e$ does not decrease distances, while it does increase at least one distance; the distance between $u$ and $v$ is 1 in $G$ and at least 2 in $G-e$. At the same time, the removal of $e$ does not increase vertices degree, while it does decrease the degree of $u$ and $v$. Similarly, adding a new edge $f=(s,t)$ to $G$ does not increase distances, while it does decrease at least one distance; the distance between $s$ and $t$ is at least 2 in $G$ and 1 in $G+f$. At the same time, the adding of $f$ does not decrease vertices degree, while it does increase the degree of $s$ and $t$.

 By the analysis above, we have the following lemma immediately which presented in \cite{hua2} for a connected graph.
\begin{lemma}
Let $G$ be a graph with $u,v \in V(G)$. If $uv \in E(G)$, then $RDD(G) > RDD(G-uv)$; If $uv \notin E(G)$, then $RDD(G) < RDD(G+uv)$.
 \end{lemma}

\section{Maximum reciprocal degree distance with given number of cut vertices}
In this section, we first introduced two edge-grafting transformations to study the mathematical properties of the reciprocal degree distance of $G$. Then using these mathematical properties, we characterize the extremal graphs with the maximum RRD-value among all the graphs of order $n$ with given cut vertices.

\begin{lemma}
Let $G_1,G_2,P_s$ be pairwise vertex-disjoint connected graphs, where $G_1$ contains an edge $uv$ such that $N_{G_1}(u)\backslash\{v\}= N_{G_1}(v)\backslash\{u\}=\{w_1,w_2,\ldots,w_k\}$ $(k\geq 1)$, $G_2$ contains a shortest path $Px_1 \ldots x_t$ from $x_1$ to $x_t$, $P_s=Pz_1z_2 \ldots z_s$, and $t\geq s+2$. Let $G$ be obtained from $G_1, G_2$ and $P_s$ by identifying $u$ with $x_1$ and $v$ with $z_1$, and let $H=G-\{z_1w_1,z_1w_2,\ldots,z_1w_k\}+\{x_2w_1,x_2w_2,\ldots,x_2w_k\}$, where $G$ and $H$ are shown in Fig. $3.1$. Then $$RDD(G)<RDD(H).$$
 \end{lemma}

 \begin{center}
\vspace{3mm}
\includegraphics[scale=.9]{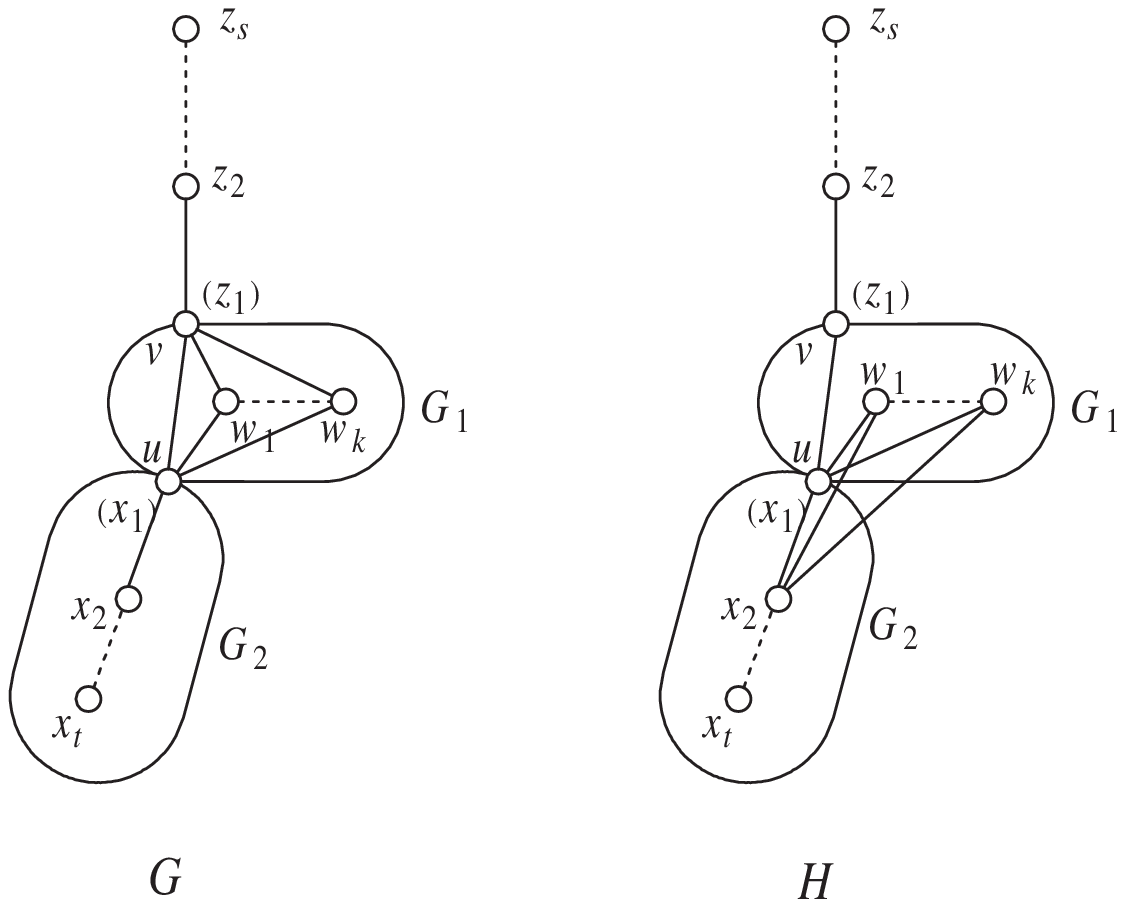}
\vspace{3mm}

{\small Fig. 3.1 \ \ The graphs $G$ and $H$ in Lemma 3.1}
\end{center}

Proof: Let $P$ be the path of $G$ obtained by connecting the path $Px_1 \ldots x_t, Puv$ and $Pz_1 \ldots z_s$, where $u=x_1$ and $v=z_1$. Partition the vertex set of $G$ as $V(G)=(V(G_1)\backslash \{u,v\})\cup (V(G_2)\backslash \{x_1,\ldots,x_t\}\cup V(P)=:S_1\cup S_2 \cup S_3$. Then from $G$ to $H$, the vertices which degrees changed only are $z_1$ and $x_2$: $d_G(z_1)=k+2$, while $d_H(z_1)=2$; and $d_G(x_2)+k=d_H(x_2)$. The vertex pairs which distance changed only are: the distance from any vertex of $S_1$ to any of $S_2$ is not increased; the distance from any vertex of $S_1$ to $z_i(i=1,2,\ldots,s)$ of $S_3$ is increased by $1$, and to $x_i(i=2,3,\ldots,t)$ is decreased by $1$.

$(1)$ Firstly, we consider the vertices of $S_1$. For any $x\in S_1$, from $G$ to $H$, the degree of $x$ is unchanged, the distance between $x$ and any other vertex of $S_1$ is unchanged, the distance between $x$ and any vertex of $S_2$ is not increased, the distance from $x$ to any of $z_i(i=1,2,\ldots,s)$ is increased by $1$, and to any of $x_i(i=2,3,\ldots,t)$ is decreased by $1$, and to the vertex $u$ is unchanged. By the analysis above and letting $d_G(x,u)=m$, we have

$$\begin{array}{lcl}
D_H(x)-D_G(x)
&\geq&(\sum_{i=0}^{t-2}\frac{1}{m+i}-\sum_{i=1}^{t-1}\frac{1}{m+i})+(\sum_{i=1}^{s}\frac{1}{m+i}-\sum_{i=0}^{s-1}\frac{1}{m+i})      \\

 & = & \frac{1}{m}-\frac{1}{m+t-1}-\frac{1}{m}+\frac{1}{m+s}\\
 & = & \frac{1}{m+s}-\frac{1}{m+t-1} \\
 & > & 0. \\
 \end{array} $$
 So, $\sum\limits_{x\in S_1}d_H(x)D_H(x)> \sum\limits_{x\in S_1} d_G(x)D_G(x).$

 $(2)$ Then we consider the vertices of $S_2$. For any $x\in S_2$, from $G$ to $H$, the degree of $x$ is unchanged, the distance between $x$ and any other vertex of $S_2$ is unchanged, the distance between $x$ and any vertex of $S_3$ is unchanged, while the distance between $x$ and any vertex of $S_1$ is not increased. By the analysis above, we have $D_H(x)-D_G(x)\geq 0.$ So, $\sum\limits_{x\in S_2}d_H(x)D_H(x)\geq \sum\limits_{x\in S_2} d_G(x)D_G(x).$

 $(3)$ Finally, we consider the vertices of $S_3$.

 From $G$ to $H$, the degree of vertex $u$ is unchanged, and the distance from $u$ to any other vertex in $G$ is unchanged. So $d_H(u)D_H(u)=d_G(u)D_G(u).$

 From $G$ to $H$, the degrees of vertices $z_2,z_3,\ldots,z_s$ are unchanged, and for any $z_i(i=2,\ldots,s)$, the distance between $z_i$ and any vertex of $S_2\cup S_3$ is unchanged, while the distance between $z_i$ and any vertex of $S_1$ is increased by $1$. Then $D_H(z_i)-D_G(z_i)=\sum\limits_{x\in S_1}\frac{1}{d_G(z_i,x)+1}-\sum\limits_{x\in S_1}\frac{1}{d_G(z_i,x)}$. Hence, $d_G(z_i)D_G(z_i)>d_H(z_i)D_H(z_i)$.

 From $G$ to $H$, the degrees of vertices $x_3,x_4,\ldots,x_t$ are unchanged, and for any $x_i(i=3,\ldots,t)$, the distance between $x_i$ and any vertex of $S_2\cup S_3$ is unchanged, while the distance between $x_i$ and any vertex of $S_1$ is decreased by $1$. Then $D_H(x_i)-D_G(x_i)=\sum\limits_{x\in S_1}\frac{1}{d_G(x_i,x)-1}-\sum\limits_{x\in S_1}\frac{1}{d_G(x_i,x)}$. Hence, $d_G(x_i)D_G(x_i)<d_H(x_i)D_H(x_i)$.

 Next we compare the change of $z_2$ and $x_3$. For any vertex $y \in S_1$, assuming $d_G(u,y)=a$, we have $d_H(y,z_2)=a+2,d_G(y,z_2)=a+1,d_H(y,x_3)=a+1,d_G(y,x_3)=a+2$. Then $$\frac{1}{d_H(y,z_2)}-\frac{1}{d_G(y,z_2)}=\frac{1}{a+2}-\frac{1}{a+1}=\frac{-1}{(a+1)(a+2)},$$
$$\frac{1}{d_H(y,x_3)}-\frac{1}{d_G(y,x_3)}=\frac{1}{a+1}-\frac{1}{a+2}=\frac{1}{(a+1)(a+2)}.$$
Notice that $d_G(z_2)=d_H(z_2)=2,d_G(x_3)=d_H(x_3)\geq 2$, we get $d_H(z_2)D_H(z_2)+d_H(x_3)D_H(x_3)\geq d_G(z_2)D_G(z_2)+d_G(x_3)D_G(x_3).$

Similarly, we can get $d_H(z_i)D_H(z_i)+d_H(x_{i+1})D_H(x_{i+1})\geq d_G(z_i)D_G(z_i)+d_G(x_{i+1})D_G(x_{i+1})$ for $i=3,\ldots,s$.

Notice that $t\geq s+2$, so $$\sum_{i=2}^{s}d_H(z_i)D_H(z_i)+\sum_{i=3}^{t}d_H(x_i)D_H(x_i)>\sum_{i=2}^{s}d_G(z_i)D_G(z_i)+\sum_{i=3}^{t}d_G(x_i)D_G(x_i).$$

In the last step, we prove $$d_H(z_1)D_H(z_1)+d_H(x_2)D_H(x_2)>d_G(z_1)D_G(z_1)+d_G(x_2)D_G(x_2).$$

Assuming $d_{G_2}(x_2)=l+2$ ($l\geq 0$), then we have
$$\begin{array}{lcl}
&&(d_H(z_1)D_H(z_1)+d_H(x_2)D_H(x_2))-(d_G(z_1)D_G(z_1)+d_G(x_2)D_G(x_2))\\
&=&2(D_H(z_1)-D_G(z_1))+(2+l)(D_H(x_2)-D_G(x_2))+k(D_H(x_2)-D_G(z_1))   \\
&=&2(\sum\limits_{x\in S_1}\frac{1}{d_G(x,z_1)+1}-\sum\limits_{x\in S_1}\frac{1}{d_G(x,z_1)})+(2+l)(\sum\limits_{x\in S_1}\frac{1}{d_G(x,x_2)-1}-\sum\limits_{x\in S_1}\frac{1}{d_G(x,x_2)})+k(D_H(x_2)-D_G(z_1))   \\
 \end{array} $$

For any $x\in S_1$, assuming $d_G(x,u)=a$, then $$\frac{1}{d_G(x,z_1)+1}-\frac{1}{d_G(x,z_1)}=\frac{1}{a+1}-\frac{1}{a}=\frac{-1}{a(a+1)},$$
$$\frac{1}{d_G(x,x_2)-1}-\frac{1}{d_G(x,x_2)}=\frac{1}{a}-\frac{1}{a+1}=\frac{1}{a(a+1)}.$$

Since $l\geq 0$, we have $$2(\sum\limits_{x\in S_1}\frac{1}{d_G(x,z_1)+1}-\sum_{x\in S_1}\frac{1}{d_G(x,z_1)})+(2+l)(\sum_{x\in S_1}\frac{1}{d_G(x,x_2)-1}-\sum_{x\in S_1}\frac{1}{d_G(x,x_2)})\geq 0$$.

For any $x\in S_1$, $d_H(x,x_2)=d_G(x,z_1)$;
For any $x\in S_2$, $d_H(x,x_2)\leq d_G(x,z_1)$, so $\frac{1}{d_H(x,x_2)}\geq \frac{1}{d_G(x,z_1)}$;

In addition,
$$\begin{array}{lcl}
\sum\limits_{y\in S_3}\frac{1}{d_H(y,x_2)} &=& 1+\frac{1}{2}+\ldots +\frac{1}{t-2}+1+\frac{1}{2}+\ldots+\frac{1}{s+1}\\
&>&1+\frac{1}{2}+\ldots+\frac{1}{s-1}+1+\frac{1}{2}+\ldots+\frac{1}{t}=\sum\limits_{y\in S_3}\frac{1}{d_G(y,z_1)}   \\
 \end{array} $$

So, we have $k(D_H(x_2)-D_G(z_1))>0.$

Thus, we proved that $$d_H(z_1)D_H(z_1)+d_H(x_2)D_H(x_2)>d_G(z_1)D_G(z_1)+d_G(x_2)D_G(x_2).$$

In view of $(1)-(3)$, we obtain $RDD(G)<RDD(H)$.
\hfill$\blacksquare$
\begin{rem}
The graphs $G$ and $H$ in Lemma $3.1$ possess the same number of cut vertices. Moreover, If taking $s=1$ in Lemma $3.1$, the edge $uv$ of $G$ will become a pendent edge of $H$.
\end{rem}

If taking $G_2=Px_1\ldots x_t$ in Lemma $3.1$, we will get the following result.
\begin{coro}
Let $G$ be a connected graph. $uv\in E(G)$ and $N_G(u)\backslash \{v\}=N_G(v)\backslash \{u\}\neq \phi$. Let $G_{s,t}$ be obtained from $G$ by attaching a path $P_t$ at $u$ and a path $P_s$ at $v$. If $t\geq s+2\geq 3$, then $RDD(G_{t,s})<RDD(G_{t-1,s+1})$.
\end{coro}

\begin{lemma}
Let $K_puK_q$ be the union of two complete graphs $K_p$ and $K_q$ sharing exactly one common vertex $u$, where $p\geq 3, q\geq 3.$ Let $G$ be obtained from $K_puK_q$ by attaching a path $P_t$ at some vertex $w_1\in V(K_p)\backslash\{u\}$ and a path $P_s$ at some vertex $v_1\in V(K_q)\backslash\{u\}$, and possibly attaching some connected graphs at other vertices of $V(K_puK_q)\backslash\{u,v_1,w_1\}$, where $t\geq s \geq 1$, and let $H$ be obtained from $G$ by deleting the edges of $K_q$ incident to $v_1$ except $v_1u$ and adding all possible edges between each of $V(K_q)\backslash \{v_1\}$ and each of $V(K_p)$, where $G$ and $H$ are shown in Fig. $3.2$. Then $RRD(G)<RRD(H)$.
\end{lemma}

 \begin{center}
\vspace{3mm}
\includegraphics[scale=.9]{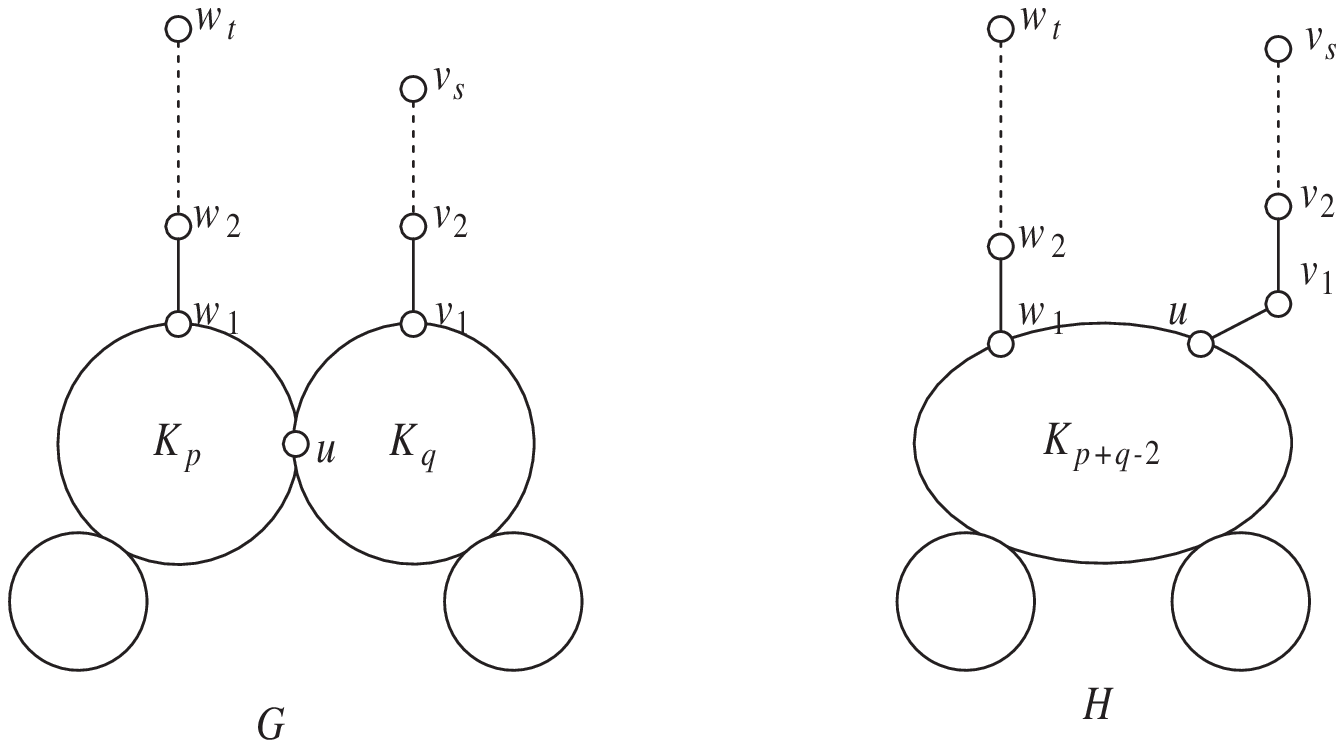}
\vspace{3mm}

{\small Fig. 3.2 \ \ The graphs $G$ and $H$ in Lemma 3.4}
\end{center}

{\it Proof:} If $t>s$, by Lemma $2.1$ and Lemma $3.1$, the result follows immediately. In the following, we discuss the case when $s=t$.

Let $S_1=\{v_1,v_2,\ldots,v_s\},S_2=\{w_1,w_2,\ldots,w_t\},S_3=\{u\},S_4=V(k_q)\backslash \{v_1,u\},S_5=V(k_p)\backslash \{w_1,u\}$, and let $S_6$ be the vertices set of the connected graphs which attached at the vertices of $V(K_q)\backslash \{u,v_1\}$, excluding the attachment points, $S_7$ be the vertices set of the connected graphs which attached at the vertices of $V(K_p)\backslash \{u,w_1\}$, excluding the attachment points. Then $V(G)$ can be partitioned as $V(G)=\cup_{i=1}^7S_i$. Observe the transformation from $G$ to $H$, the degree of $v_1$ changes from $q$ to $2$, the degree of $w_1$ changes from $p$ to $p+q-2$, the degree of any vertex in $S_4$ is increased by $p-2$, the degree of any vertex in $S_5$ is increased by $q-2$, while the degrees of any other vertex is unchanged; the distance between any vertex of $S_4\cup S_6$ and any of $S_1$ are increased by $1$, the distance between any vertex of $S_4\cup S_6$ and any of $S_2\cup S_5\cup S_7$ is decreased by $1$, while the distance between any other two vertices is not changed.

(1) Firstly, we consider the vertices $v_2$ and $w_2$.
$$\begin{array}{lcl}
&&d_H(v_2)D_H(v_2)-d_G(v_2)D_G(v_2)+d_H(w_2)D_H(w_2)-d_G(w_2)D_G(w_2)\\
&=&2(\sum_{x\in S_4\cup S_6}\frac{1}{d_H(v_2,x)}-\sum_{x\in S_4\cup S_6}\frac{1}{d_G(v_2,x)})+2(\sum_{x\in S_4\cup S_6}\frac{1}{d_H(w_2,x)}-\sum_{x\in S_4\cup S_6}\frac{1}{d_G(w_2,x)}) \\
&=&2(\sum_{x\in S_4\cup S_6}\frac{1}{d_G(v_2,x)+1}-\sum_{x\in S_4\cup S_6}\frac{1}{d_G(v_2,x)})+2(\sum_{x\in S_4\cup S_6}\frac{1}{d_G(w_2,x)-1}-\sum_{x\in S_4\cup S_6}\frac{1}{d_G(w_2,x)}) \\
 \end{array} $$
For any vertex $x\in S_4\cup S_6$, let $d_G(w_2,x)=a$, then $d_G(v_2,x)=a-1$. So,
$$(\frac{1}{d_G(v_2,x)+1}-\frac{1}{d_G(v_2,x)})+(\frac{1}{d_G(w_2,x)-1}-\frac{1}{d_G(w_2,x)})=\frac{1}{a}-\frac{1}{a-1}+\frac{1}{a-1}-\frac{1}{a}=0.$$
Hence, $d_H(v_2)D_H(v_2)+d_H(w_2)D_H(w_2)-d_G(v_2)D_G(v_2)-d_G(w_2)D_G(w_2)=0.$
Similarly, for any $v_i$ and $w_i$,$i=3,\ldots,s$, $d_G(v_i)D_G(v_i)+d_G(w_i)D_G(w_i)=d_H(v_i)D_H(v_i)+d_H(w_i)D_H(w_i).$

(2) For the vertex $u$, $d_G(u)=d_H(u)$, and the distance from $u$ to any other vertex is unchanged, so $d_G(u)D_G(u)=d_H(u)D_H(u).$

(3) For any vertex $x\in S_4$, the degree of $x$ is increased by $p-2$, the distance from $x$ to any vertex of $S_1$ is increased by $1$, the distance from $x$ to any vertex of $S_2 \cup S_5 \cup S_7$ is decreased by $1$. For any $i=1,2,\ldots,s$, let $d_G(x,v_i)=a_i$, then $d_G(x,w_i)=a+1$, so $\frac{1}{d_H(x,v_i)}+\frac{1}{d_H(x,w_i)}-\frac{1}{d_G(x,v_i)}-\frac{1}{d_G(x,w_i)}=\frac{1}{a_i+1}+\frac{1}{a_i}-\frac{1}{a_i}-\frac{1}{a_i+1}=0.$ Hence, $d_G(x)D_G(x)<d_H(x)d_H(x)$.

(4) For any vertex $x\in S_5$, the degree of $x$ is increased by $q-2$, the distance from $x$ to any vertex of $S_4\cup S_6$ is decreased by $1$, and the distance from $x$ to any other vertex is unchanged. Hence, $d_G(x)D_G(x)<d_H(x)D_H(x).$

(5) For any vertex $x\in S_6$, the degree of $x$ is unchanged, the distance from $x$ to any vertex of $S_1$ is increased by $1$, the distance from $x$ to any vertex of $S_2 \cup S_5 \cup S_7$ is decreased by $1$. By similar discussion to the vertex of $S_4$, we can get $d_G(x)D_G(x)<d_H(x)d_H(x)$.

(6) For any vertex $x\in S_7$, the degree of $x$ is unchanged, the distance from $x$ to any vertex of $S_4\cup S_6$ is decreased by $1$, and the distance from $x$ to any other vertex is unchanged. Hence, $d_G(x)D_G(x)<d_H(x)D_H(x).$

(7) In the last step, we concentrate on the vertices $v_1$ and $w_1$. From $G$ to $H$, the degree of $v_1$ is changed from $q$ to $2$, the degree of $w_1$ is changed from $p$ to $p+q-2$, the distance from $v_1$ to any vertex of $S_4\cup S_6$ is increased by $1$, the distance from $v_1$ to any other vertex is unchanged, the distance from $w_1$ to any vertex of $S_4\cup S_6$ is decreased by $1$, the distance from $w_1$ to any other vertex is unchanged. For simplicity, let $A=S_4\cup S_6, B=V(G)-S_4\cup S_6$. Thus, we have
$$\begin{array}{lcl}
&&d_H(v_1)D_H(v_1)-d_G(v_1)D_G(v_1)+d_H(w_1)D_H(w_1)-d_G(w_1)D_G(w_1)\\
&=&2(\sum\limits_{x\in A}\frac{1}{d_H(v_1,x)}+\sum\limits_{x\in B-\{v_1\}}\frac{1}{d_H(v_1,x)})-q(\sum\limits_{x\in A}\frac{1}{d_G(v_1,x)}+\sum\limits_{x\in B-\{v_1\}}\frac{1}{d_G(v_1,x)})\\
&&+(p+q-2)(\sum\limits_{x\in A}\frac{1}{d_H(w_1,x)}+\sum\limits_{x\in B-\{w_1\}}\frac{1}{d_H(w_1,x)})-p(\sum\limits_{x\in A}\frac{1}{d_G(w_1,x)}+\sum\limits_{x\in B-\{w_1\}}\frac{1}{d_G(w_1,x)})\\
&=&2\sum\limits_{x\in A}\frac{1}{d_H(v_1,x)}-q\sum\limits_{x\in A}\frac{1}{d_G(v_1,x)}+(p+q-2)\sum\limits_{x\in A}\frac{1}{d_H(w_1,x)}-p\sum\limits_{x\in A}\frac{1}{d_G(w_1,x)}\\
&&+(2-q)\sum\limits_{x\in B-\{v_1\}}\frac{1}{d_G(v_1,x)}+(q-2)\sum\limits_{x\in B-\{w_1\}}\frac{1}{d_G(w_1,x)}\\
 \end{array} $$
 For any $x\in A$, let $d_G(v_1,x)=a$, then $d_G(w_1,x)=a+1,d_H(v_1,x)=a+1,d_H(w_1,x)=a$, thus,
 $$\begin{array}{lcl}
2\frac{1}{d_H(v_1,x)}-q\frac{1}{d_G(v_1,x)}+(p+q-2)\frac{1}{d_H(w_1,x)}-p\frac{1}{d_G(w_1,x)}&=&2\frac{1}{a+1}-q\frac{1}{a}+(p+q-2)\frac{1}{a}-p\frac{1}{a+1}\\
&=&(2-p)\frac{1}{a+1}+(p-2)\frac{1}{a}\\
&=&(p-2)(\frac{1}{a}-\frac{1}{a+1})\\
&>&0.\\
 \end{array} $$
Hence, $$2\sum\limits_{x\in A}\frac{1}{d_H(v_1,x)}-q\sum\limits_{x\in A}\frac{1}{d_G(v_1,x)}+(p+q-2)\sum\limits_{x\in A}\frac{1}{d_H(w_1,x)}-p\sum\limits_{x\in A}\frac{1}{d_G(w_1,x)}>0.$$
In addition, for any vertex pairs $v_i$ and $w_i$, $(i=2,3,\ldots,s)$, $\frac{1}{d_G(w_1,v_i)}+\frac{1}{d_G(w_1,w_i)}=\frac{1}{d_G(v_1,v_i)}+\frac{1}{d_G(v_1,w_i)}$; for the vertex $u$, $d_G(w_1,u)=d_G(v_1,u)$; while for any vertex $x\in S_5\cup S_7$, $d_G(w_1,x)<d_G(v_1,x)$. Thus
$$\begin{array}{lcl}
&&(2-q)\sum\limits_{x\in B-\{v_1\}}\frac{1}{d_G(v_1,x)}+(q-2)\sum\limits_{x\in B-\{w_1\}}\frac{1}{d_G(w_1,x)}\\
&=&(2-q)(\sum\limits_{x\in B-\{v_1,w_1\}}\frac{1}{d_G(v_1,x)}+\frac{1}{d_G(v_1,w_1)})+(q-2)(\sum\limits_{x\in B-\{w_1,v_1\}}\frac{1}{d_G(w_1,x)}+\frac{1}{d_G(w_1,v_1)})\\
&=&(q-2)(\sum\limits_{x\in B-\{w_1,v_1\}}\frac{1}{d_G(w_1,x)}-\sum\limits_{x\in B-\{w_1,v_1\}}\frac{1}{d_G(v_1,x)})\\
&=&(q-2)(\sum\limits_{x\in S_5\cup S_7}\frac{1}{d_G(w_1,x)}-\sum\limits_{x\in B_5 \cup B_7}\frac{1}{d_G(v_1,x)})\\
&>&0.
 \end{array} $$
Therefore, $$d_H(v_1)D_H(v_1)-d_G(v_1)D_G(v_1)+d_H(w_1)D_H(w_1)-d_G(w_1)D_G(w_1)>0$$.

Combining $(1)-(7)$, the result follows.
\hfill$\blacksquare$

\begin{rem} The graphs $G$ and $H$ in Lemma $3.4$ possess the same number of cut vertices. Moreover, if taking $s=1$, the edge $uv_1$ of $G$ becomes a pendent edge of $H$.
\end{rem}

\begin{theorem} For any $G\in \mathcal{G}_{n,k}$,  where $0\leq k \leq n-2$, $$RDD(G)\leq RDD(G_{n,k}),$$ with equality holds if and only if $G\cong G_{n,k}$.
\end{theorem}
{\it Proof:} Let $G_0$ be a graph with the maximal reciprocal degree distance among all the graphs with $n$ vertices and $k$ cut vertices. If $k=0$, then by Lemma 2.1, $G_0\cong K_n \cong G_{n,0}$. Suppose in the following that $1\leq k \leq n-2$.

Claim 1: $G_0$ is connected.

 If $G_0$ is disconnected, then $G_0$ has at least two components. Let $z$ be a cut vertex of $G_0$. Then $z$ is also a cut vertex of some component, say $H_1$, of $G_0$. Let $H_2$ be another component of $G_0$. If there is a cut vertex, say $z'$, in $H_2$, then $G_0+{zz'}$ possesses $k$ cut vertices, and by Lemma 2.1, $RDD(G_0)<RDD(G_0+{zz'})$, a contradiction. If there is no cut vertex in $H_2$, then denote by $G_0'$ the graph obtained from $G_0$ by adding the edges between $z$ and all vertices of $H_2$. Thus $G_0'$ also possesses $k$ cut vertices, and by Lemma 2.1, $RDD(G_0)<RDD(G_0')$, a contradiction again. Hence $G_0$ is connected.

By Lemma 2.1, each block of $G_0$ is complete, and each cut vertex of $G_0$ is contained exactly in two blocks.
If each block of $G_0$ has exactly two vertices, i.e., each block is a single edge, then $G_0$ is a tree with maximum degree two, i.e., $G_0\cong P_n\cong G_{n,n-2}$. Suppose in the following that there is at least one block of $G_0$ with at least three vertices.

Claim 2: If $G_0\neq G_{n,1}$, then each pendent block of $G_0$ is an edge.

If $B_1$ is a pendent block of $G_0$ and $|V(B_1)|>2$, we assume $u$ is a vertex different from the unique cut vertex, say $w$, of $B_1$. Denote by $B_2$ the block adjacent to $B_1$. By deleting the edges between $u$ and $V(B_1)-\{u,w\}$, and adding all the edges between $V(B_1)-\{u,w\}$ and $V(B_2)-w$, we obtain a new graph $G_0'$. Notice that the number of cut vertices of $G_0'$ is also $k$, and by remark $3.2$ (if $|V(B_2)|=2$) and remark $3.5$ (if $|V(B_2)|>2$), we have $RDD(G_0)<RDD(G_0')$, a contradiction.

Choose a pendent path, say $P_s$ at $v$, with minimal length in $G_0$. Obviously, $v$ lies in some block, say $B$, of $G_0$ with at least three vertices. Note that $v$ is not a cut vertex of $G_0$ if $s=1$.

Claim 3: The component attached at any vertex of $B$ is a path(possibly being trivial).

For $x\in V(B)$, let $H^{(x)}$ be the component of $G-E(B)$ containing $x$. Obviously, $H^{(v)}\cong P_s$. Suppose $u$ is an arbitrary  vertex of $B$ and $u \neq v$. Obviously, $N_B(v)\backslash \{u\}=N_B(u)\backslash \{v\}$. Let $G^*$ be the component of $G-((E(H^{(u)})\cup E(P_s))$ containing $u$, which surely contains the block $B$.

Suppose that $H^{(u)}$ is not a (possibly trivial) path. Then $H^{(u)}$ contains a block with at least three vertices. By the proof of Claim 2, $H^{(u)}$ must contain a nontrivial pendant path $P_t$ attached at some nontrivial block $B_0$ of $H^{u}$, where $s\geq t$. Therefore $H^{u}$ contains a shortest path $P_r$ from $u$ to the pendent vertex of $P_t$, where $r\geq t+1\geq s+1$. If $s=1$, then by Remark $3.2$, we may get another graph with $n$ vertex and $k$ cut vertices, which has a larger reciprocal degree distance, a contradiction. If $s>1$ and $r\geq s+2$, then by Lemma $3.1$, we also get a contradiction. So in the following we only need to consider the case:$s>1$ and $r=s+1$. In this case, $B_0$ share with $B$ the common vertex $u$, and $H^{(u})$ is obtained from $B_0$ by attaching $P_s$ at each of its vertices except $u$. Applying Lemma $3.4$, we can get another graph of order $n$ with $k$ cut vertices, which has a larger reciprocal degree distance, a contradiction. Therefore $H^{(u)}$ is a pendent path attached at $u$ which contains at least $s$ vertices.

Claim 4: All paths attached at the vertices of $B$ have almost equal lengths.

Obviously, $t\geq s$. If $t \geq s+2$, then by Corollary $3.3$, we may get another graph with $n$ vertices and $k$ cut vertices, which has a larger reciprocal degree distance, a contradiction. So $H^{(u)}\cong P_s$ or $P_{s+1}$.

To sum up, we get $G\cong \mathbf{G}_{n,k}$.
\hfill$\blacksquare$

\section{Maximum reciprocal degree distance with given number of cut edges}
Similar to section 3, we first introduced two edge-grafting transformations to study the mathematical properties of the reciprocal degree distance of $G$. Then using these mathematical properties, we characterize the extremal graphs with the maximum RRD-value among all the graphs of order $n$ with given cut edges. In addition, we also provide an upper bounds on the reciprocal degree distance in terms of the number of cut edges. The following lemma is a special case of Theorem $2.1$ in \cite{li1}.
\begin{lemma}
Let $w_1w_2\in E(G)$ be a cut edge in $G$, and $G-w_1w_2=G_1\cup G_2$ where $G_i$ is nontrivial and $w_i\in V(G_i)$ for $i=1,2$. Assume that $H$ is a graph obtained from $G$ by identifying $w_1$ with $w_2$ (the new vertex is labeled as $w$) and attaching at $w$ a pendent vertex $w_0$. $G$ and $H$ are shown in Fig. $4.1$. Then $RDD(G)<RDD(H)$.
\end{lemma}

\begin{center}
\vspace{3mm}
\includegraphics[scale=.9]{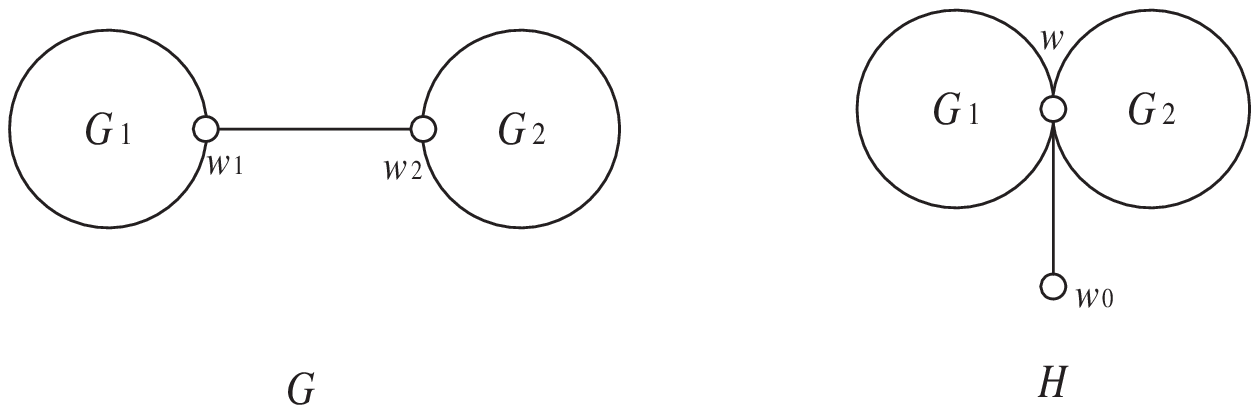}
\vspace{3mm}

{\small Fig. 4.1 \ \ The graphs $G$ and $H$ in Lemma 4.1}
\end{center}

\begin{lemma}
Let $G_0,G_1,G_2$ be pairwise vertex-disjoint connected graphs and $u,v \in V(G_0)$ such that $N_{G_0}(u)\backslash \{v\}=N_{G_0}(v)\backslash \{u\}$, $w_1\in V(G_1)$, $w_2\in V(G_2)$. Let $H$ be the graph obtained from $G_0,G_1,G_2$ by identifying $u$ with $w_1$ and $v$ with $w_2$, respectively. Let $H_1$ be the graph obtained from $G_0,G_1,G_2$ by identifying three vertices $u,w_1,w_2$, and let $H_2$ be the graph obtained from $G_0,G_1,G_2$ by identifying three vertices $v,w_1,w_2$. $H,H_1$ and $H_2$ are shown in Fig. $4.2$. Then we have $RDD(H_i)> RDD(H)$ for $i=1,2$.
\end{lemma}

\begin{center}
\vspace{3mm}
\includegraphics[scale=.9]{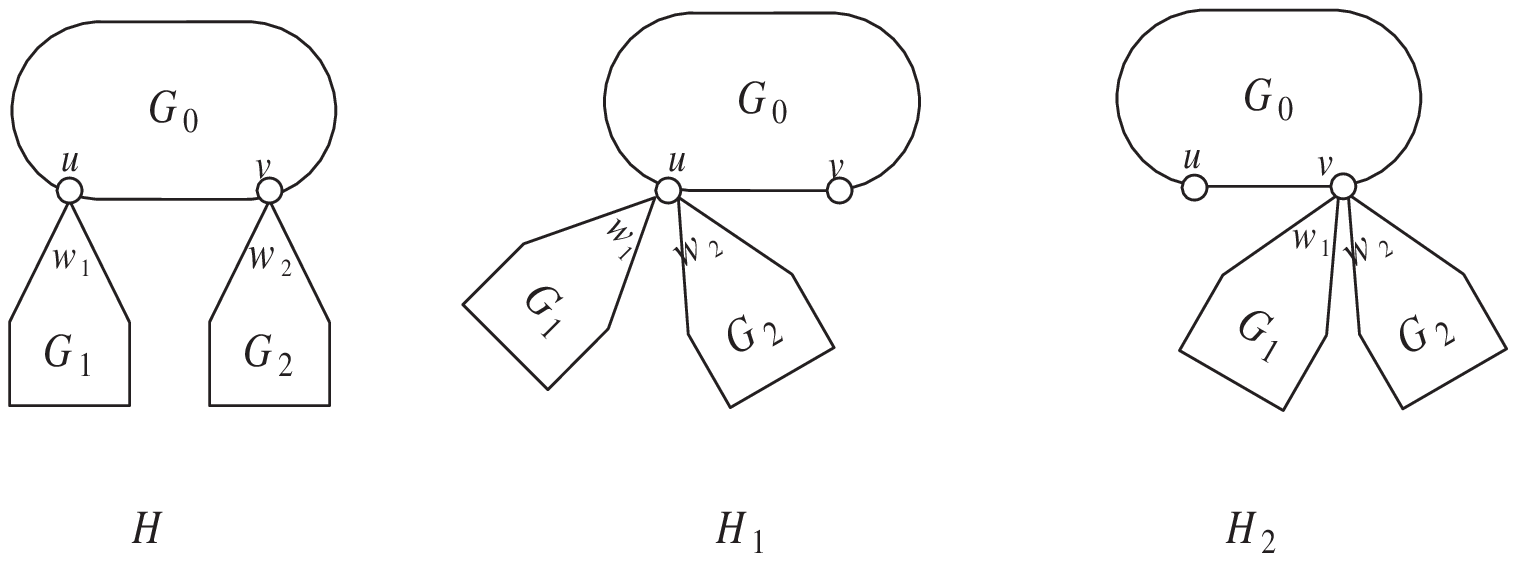}
\vspace{3mm}

{\small Fig. 4.2 \ \ The graphs $H,H_1$ and $H_2$ in Lemma 4.2}
\end{center}

{\it Proof:} For simplicity, we denote $G_0-\{u,v\}$ by $G_0^*$, $G_1-w_1$ by $G_1^*$ and $G_2-w_2$ by $G_2^*$. Then $V(H)=V(H_1)=V(H_2)=V(G_0^*)\cup V(G_1^*)\cup V(G_2^*)\cup \{u,v\}$. Obviously, $V(G_0^*),V(G_1^*),V(G_2^*)$ and$\{u,v\}$ are four vertex sets disjoint in pair. Since $N_{G_0}(u)\backslash \{v\}=N_{G_0}(v)\backslash \{u\}$, we have $d_{G_0}(u)=d_{G_0}(v)$ and for any $x\in V(G_0^*)$, $d_{G_0}(u,x)=d_{G_0}(v,x)$. Note that from $H$ to $H_i(i=1,2)$, the vertices which degree changed only are $u$ and $v$. Since in $H, H_1$ or $H_2$, $u$ and $v$ have the same distance. For simplicity, we denote by $d(u,v)$ the distance between $u$ and $v$ in $H, H_1$ or $H_2$. Similarly, $d(x,u)$ ($d(x,v)$, respectively) denotes the distance between $x$ and $u$ ($v$, respectively) for any $x\in V(G_0^*)$, $d(w_1,y)$ denotes  the distance between $w_1$ and $y$ for any $y\in V(G_1^*)$, and $d(w_2,z)$ denotes  the distance between $w_2$ and $z$ for any $z\in V(G_2^*)$. Therefore,

$$\begin{array}{lcl}
RDD(H_1)-RDD(H)\\
=\sum\limits_{x\in V(G_0^*)}[d(x)\sum\limits_{z\in V(G_2^*)}(\frac{1}{d_{H_1}(x,z)}-\frac{1}{d_{H}(x,z)})]+\sum\limits_{y\in V(G_1^*)}[d(y)\sum\limits_{z\in V(G_2^*)}(\frac{1}{d_{H_1}(y,z)}-\frac{1}{d_{H}(y,z)})]\\
+\sum\limits_{z\in V(G_2^*)}[d(z)(\sum\limits_{x\in V(G_0^*)}(\frac{1}{d_{H_1}(x,z)}-\frac{1}{d_{H}(x,z)})+\sum\limits_{y\in V(G_1^*)}(\frac{1}{d_{H_1}(y,z)}-\frac{1}{d_{H}(y,z)})+\frac{1}{d_{H_1}(z,u)}-\frac{1}{d_{H}(z,u)}+\frac{1}{d_{H_1}(z,v)}-\frac{1}{d_{H}(z,v)})]\\
+d_{H_1}(u)(\sum\limits_{x\in V(G_0^*)}\frac{1}{d(x,u)}+\sum\limits_{y\in V(G_1^*)}\frac{1}{d(w_1,y)}+\sum\limits_{z\in V(G_2^*)}\frac{1}{d(w_2,z)}+\frac{1}{d(u,v)})\\
-d_{H}(u)(\sum\limits_{x\in V(G_0^*)}\frac{1}{d(x,u)}+\sum\limits_{y\in V(G_1^*)}\frac{1}{d(w_1,y)}+\sum\limits_{z\in V(G_2^*)}\frac{1}{d(u,v)+d(w_2,z)}+\frac{1}{d(u,v)})\\
+d_{H_1}(v)(\sum\limits_{x\in V(G_0^*)}\frac{1}{d(x,v)}+\sum\limits_{y\in V(G_1^*)}\frac{1}{d(w_1,y)+d(u,v)}+\sum\limits_{z\in V(G_2^*)}\frac{1}{d(w_2,z)+d(u,v)}+\frac{1}{d(u,v)})\\
-d_{H}(v)(\sum\limits_{x\in V(G_0^*)}\frac{1}{d(x,v)}+\sum\limits_{y\in V(G_1^*)}\frac{1}{d(w_1,y)+d(u,v)}+\sum\limits_{z\in V(G_2^*)}\frac{1}{d(w_2,z)}+\frac{1}{d(u,v)})\\
>(d_{G_0}(u)+d_{G_1}(w_1)-d_{G_0}(v))\sum\limits_{z\in V(G_2^*)}(\frac{1}{d(w_2,z)}-\frac{1}{d(u,v)+d(w_2,z)})
+d_{G_2}(w_2)\sum\limits_{x\in V(G_0^*)}(\frac{1}{d(x,u)}-\frac{1}{d(x,v)})\\
+\sum\limits_{x\in V(G_0^*)}\sum\limits_{z\in V(G_2^*)}(d(x)+d(z))(\frac{1}{d(x,u)+d(w_2,z)}-\frac{1}{d(x,v)+d(w_2,z)})\\
=d_{G_1}(w_1)\sum\limits_{z\in V(G_2^*)}(\frac{1}{d(w_2,z)}-\frac{1}{d(u,v)+d(w_2,z)})\\
>0\\
\end{array} $$

Similarly, we have
$$\begin{array}{lcl}
RDD(H_2)-RDD(H)\\
>(d_{G_0}(v)+d_{G_2}(w_2)-d_{G_0}(u))\sum\limits_{y\in V(G_1^*)}(\frac{1}{d(w_1,y)}-\frac{1}{d(u,v)+d(w_1,y)})
+d_{G_1}(w_1)\sum\limits_{x\in V(G_0^*)}(\frac{1}{d(x,v)}-\frac{1}{d(x,u)})\\
+\sum\limits_{x\in V(G_0^*)}\sum\limits_{y\in V(G_1^*)}(d(x)+d(y))(\frac{1}{d(x,v)+d(w_1,y)}-\frac{1}{d(x,u)+d(w_1,y)})\\
=d_{G_2}(w_2)\sum\limits_{y\in V(G_1^*)}(\frac{1}{d(w_1,y)}-\frac{1}{d(u,v)+d(w_1,y)})\\
>0\\
\end{array} $$
The result follows.
\hfill$\blacksquare$

\begin{theorem}
 For any $G\in \overline{\mathcal{G}}_{n,k}$, $$RDD(G)\leq n^3-(\frac{5}{2}k+2)n^2+(2k^2+\frac{11}{2}k+1)n-(\frac{1}{2}k^3+2k^2+\frac{5}{2}k),$$ with equality holds if and only if $G\cong \overline{G}_{n,k}$.
\end{theorem}
{\it Proof:} Let $G_0$ be a graph with the minimum reciprocal degree distance in $\overline{\mathcal{G}}_{n,k}$ and $E_1=\{e_1,e_2,\ldots,e_k\}$ be the set of cut edges of $G_0$. Firstly, by Lemma $2.1$, we can get each component of $G_0-E_1$ is a clique. In addition, by Lemma $4.1$, $e_1,e_2,\ldots,e_k$ must be the pendent edges in $G_0$. Hence, $G_0$ must be the graph obtain from $K_{n-k}$ by attaching $k$ pendent edges to some vertices. Finally, by Lemma $4.2$, all these pendent edges in $G_0$ must be attached to one common vertex. Thus $G_0\cong \overline{G}_{n,k}$.

In the following we only need to calculate $RRD(\overline{G}_{n,k})$.
By the structure of $\overline{G}_{n,k}$, we can get
$$\begin{array}{lcl}
RDD(\overline{G}_{n,k})&=&k(1+\frac{n-2}{2})+(n-1)^2+(n-k-1)[(n-k-1)(n-k-1+\frac{k}{2})]\\
&=&n^3-(\frac{5}{2}k+2)n^2+(2k^2+\frac{11}{2}k+1)n-(\frac{1}{2}k^3+2k^2+\frac{5}{2}k).\\
\end{array} $$

This completes the proof.
\hfill$\blacksquare$

\end{document}